\documentclass[11pt]{article}
\usepackage{amsmath,amssymb,amsfonts}
\usepackage[latin1]{inputenc}
\usepackage{epsfig,psfrag}
\numberwithin{equation}{section}
\newtheorem{thm}{Theorem}[section]

\newtheorem{alem}[thm]{Lemma}
\newtheorem{aprop}[thm]{Proposition}

\newtheorem{arem}[thm]{Remark}
\newenvironment{adem}[1][]%
   {\ \\ {\bf Proof #1. }}%
   {\hfill\mbox{\rule{2 true mm}{3 true mm}}}
   {\ \\ {\bf Example #1. }}%
   {\hfill\mbox{\rule{2 true mm}{3 true mm}}}
\newcommand{\R}{\mathbb{R}}

\newcommand{\E}{{\mathbb E}}

\newcommand{\mycomments}[1]{  }

\newcommand{\be}{ \begin{eqnarray*}  }
\newcommand{\ee}{ \end{eqnarray*}  }

 \title{A remark on the optimal transport between two probability measures sharing the same copula}
\author{A. Alfonsi, B. Jourdain\thanks{Universit\'e Paris-Est, CERMICS, Projet MathFi
    ENPC-INRIA-UMLV, 6 et 8 avenue Blaise Pascal, 77455 Marne La Vall\'ee, Cedex
    2, France, e-mails : alfonsi@cermics.enpc.fr, jourdain@cermics.enpc.fr. This research benefited
    from the support of the ``Chaire Risques Financiers'', Fondation du
    Risque and the French National Research Agency (ANR) under the program
 ANR-12-BLAN Stab.
}}

\begin{document}
\maketitle
\begin{abstract}
We are interested in the Wasserstein distance between two probability measures on $\R^n$ sharing the same copula $C$. The image of the probability measure $dC$  by the vectors of pseudo-inverses of marginal distributions is a natural generalization of the coupling known to be optimal in dimension $n=1$. It turns out that for cost functions $c(x,y)$ equal to the $p$-th power of the $L^q$ norm of $x-y$ in $\R^n$, this coupling is optimal only when $p=q$ i.e. when $c(x,y)$ may be decomposed as the sum of coordinate-wise costs.
\end{abstract}

\noindent {\bf Keywords:} Optimal transport, Copula, Wasserstein distance, Inversion of the cumulative distribution function.

\noindent {\bf AMS Classification (2010):} 60E05, 60E15.

\section{Optimal transport between two probability measures sharing the same copula}
Given two probability measures $\mu$ and $\rho$, the optimal transport theory aims at minimizing $\int c(x,y)\nu(dx,dy)$ over all couplings $\nu$ with first marginal $\nu\circ((x,y)\mapsto x)^{-1}=\mu$ and second marginal $\nu\circ((x,y)\mapsto y)^{-1}=\rho$ for a measurable non-negative cost function $c$. We use the notation $\nu\hspace{-1mm}<^\mu_\rho$ for such couplings. In the present note, we are interested in the particular case of the so-called Wasserstein distance between two probability measures $\mu$ and $\rho$ on $\R^n$ :
\begin{equation}\label{wass}
{\cal W}_{p,q}(\mu,\rho)=\inf_{\nu<^\mu_\rho}\left(\int_{\R^n\times\R^n}\|x-y\|_q^p\nu(dx,dy)\right)^{1/p}
\end{equation}
obtained for the choice $c(x,y)=\|x-y\|^p_q$. Here $\R^n$ is endowed with the norm $\|(x_1,\hdots,x_n)\|_q=\left(\sum_{i=1}^n|x_i|^q\right)^{1/q}$ for $q\in[1,+\infty)$ whereas $p\in[1,+\infty)$ is the power of this norm in the cost function.\par
In dimension $n=1$, $\|x\|_q=|x|$ so that the Wasserstein distance does not depend on $q$ and is simply denoted by ${\cal W}_p$. Moreover, the optimal transport is given by the inversion of the cumulative distribution functions : whatever $p\in[1,+\infty)$, the optimal coupling is the image of the Lebesgue measure on $(0,1)$ by $u\mapsto (F^{-1}_\mu(u),F^{-1}_\rho(u))$ where for $u\in (0,1)$, $F^{-1}_\mu(u)=\inf\{x\in\R:\mu((-\infty,x])\geq u\}$ and $F^{-1}_\rho(u)=\inf\{x\in\R:\rho((-\infty,x])\geq u\}$ (see for instance Theorem 3.1.2 in~\cite{raru}). This implies that ${\cal W}^p_{p}(\mu,\rho)=\int_{(0,1)}|F^{-1}_\mu(u)-F^{-1}_\rho(u)|^pdu$.

In higher dimensions, according to Sklar's theorem (see for instance Theorem 2.10.11 in Nelsen~\cite{nelsen}), 
$$\mu\left(\prod_{i=1}^n (-\infty,x_i]\right)=C\left(\mu_1((-\infty,x_1]),\hdots,\mu_n((-\infty,x_n])\right)$$
where we denote by $\mu_i=\mu\circ((x_1,\hdots,x_n)\mapsto x_i)^{-1}$ the $i$-th marginal of $\mu$ and $C$ is a copula function i.e. $C(u_1,\hdots,u_n)=m\left(\prod_{i=1}^n [0,u_i]\right)$ for some probability measure $m$ on $[0,1]^n$ with all marginals equal to the Lebesgue measure on $[0,1]$. The copula function $C$ is uniquely determined on the product of the ranges of the marginal cumulative distribution functions $x_i\mapsto\mu_i((-\infty,x_i])$. In particular, when the marginals $\mu_i$ do not weight points, the copula~$C$ is uniquely determined. Sklar's theorem shows that the dependence structure associated with $\mu$ is encoded in the copula function $C$. Last, we give the well-known Fr\'echet-Hoeffding bounds 
$$\forall u_1,\dots,u_n \in [0,1], \ C^-_n(u_1,\dots,u_n)\le C(u_1,\dots,u_n)\le C^+_n(u_1,\dots,u_n)$$
that hold for any copula function $C$ with $C^+_n(u_1,\dots,u_n)=\min(u_1,\dots,u_n)$ and $C^-_n(u_1,\dots,u_n)=(u_1+\dots+u_n-n+1)^+$ (see Nelsen~\cite{nelsen}, Theorem 2.10.12 or Rachev and R\"uschendorf~\cite{raru}, section 3.6). We recall that the copula $C^+_n$ is the $n$-dimensional cumulative distribution function of the image of the Lebesgue measure on $[0,1]$ by $\R\ni x\mapsto (x,\dots,x)\in\R^n$. Also the copula $C^-_2$ is the $2$-dimensional cumulative distribution function of the image of the Lebesgue measure on $[0,1]$ by $\R\ni x\mapsto (x,1-x)\in\R^2$ and, for $n\ge 3$, $C^-_n$ is not a copula.

 \par
In dimension $n=1$, the unique copula function is $C(u)=u$ and therefore the optimal coupling between $\mu$ and $\rho$, which necessarily share this copula, is the image of the probability measure $dC$ by $u\mapsto (F^{-1}_\mu(u),F^{-1}_\rho(u))$. It is therefore natural to wonder whether, when $\mu$ and $\rho$ share the same copula~$C$ in higher dimensions, the optimal coupling is still the image of the probability measure $dC$ by $(u_1,\hdots,u_n)\mapsto (F_{\mu_1}^{-1}(u_1),\hdots,F_{\mu_n}^{-1}(u_n),F_{\rho_1}^{-1}(u_1),\hdots,F_{\rho_n}^{-1}(u_n))$. We denote by $\mu \diamond \rho$ this probability law on~$\R^{2n}$. It turns out that the picture is more complicated than in dimension one because of the choice of the index $q$ of the norm.
\begin{aprop}\label{mr}Let $n\geq 2$, $\mu$ and $\rho$ be two probability measures on $\R^n$ sharing the same copula~$C$ and ${\cal W}_{p,q}(\mu,\rho)=\inf_{\nu<^\mu_\rho}\left(\int_{\R^n\times\R^n}\|x-y\|_q^p\nu(dx,dy)\right)^{1/p}$.
\begin{itemize}
   \item If $p=q$, then an optimal coupling between $\mu$ and $\rho$ is given by $\nu=\mu \diamond \rho$ 
     and
$${\cal W}^p_{p,p}(\mu,\rho)=\int_{[0,1]^n}\sum_{i=1}^n|F_{\mu_i}^{-1}(u_i)-F_{\rho_i}^{-1}(u_i)|^p dC(u_1,\hdots,u_n)=\int_{[0,1]}\sum_{i=1}^n|F_{\mu_i}^{-1}(u)-F_{\rho_i}^{-1}(u)|^pdu.$$
\item If $p\neq q$, the coupling $\mu \diamond \rho$ is in general no longer optimal. For
$p<q$, if $C \neq C_n^+$, we can construct probability measures~$\mu$ and $\rho$ on $\R^n$ admitting~$C$ as their unique copula such that
$$\left(\int_{\R^n\times\R^n}\|x-y\|_q^p \mu \diamond \rho(dx,dy)\right)^{1/p} >  {\cal W}_{p,q}(\mu,\rho). $$
For $p>q$, the same conclusion holds if $n\ge 3$ or $n=2$ and $C \neq C_2^-$.
\end{itemize}
\end{aprop}

\begin{arem}\label{rk_puccetti_Scarsini}
 Let $\mu$ and $\rho$ be two probability measures on $\R^n$ and $\nu\hspace{-1mm}<^\mu_\rho$. For $n=1$, $\nu$ is said to be comonotonic if $\nu((-\infty,x],(-\infty,y])=C^+_2(\mu((-\infty,x]),\rho((-\infty,y]))$.  Puccetti and Scarsini~\cite{pusc} investigate several extensions of this notion for $n\ge 2$. In particular, they say that $\nu$ is $\pi$-comonotonic (resp. $c$-comonotonic) if $\mu$ and $\rho$ have a common copula and $\nu=\mu \diamond \rho$ (resp. $\nu$ maximizes $\int_{\R^n\times\R^n} \langle x, y \rangle \tilde{\nu}(dx,dy)$ over all the coupling measures $\tilde{\nu}\hspace{-1mm}<^\mu_\rho$). Looking at some connections between their different definitions of comonotonicity, they show in Lemma 4.4 that $\pi$-comonotonicity implies $c$-comonotonicity. Since 
$$\int_{\R^n\times\R^n}\|x-y\|_2^2\tilde{\nu}(dx,dy) =\int_{\R^n}\|x\|_2^2\mu(dx)+\int_{\R^n}\|y\|_2^2\rho(dy)  -2\int_{\R^n\times\R^n} \langle x, y \rangle \tilde{\nu}(dx,dy),$$
this yields our result in the case $p=q=2$.
\end{arem}

\section{Proof of Proposition \ref{mr}}
The optimality in the case $q=p$, follows by choosing $d_1=\hdots=d_n=d'_1=\hdots=d'_n=d''_1=\hdots=d''_n=1$, $c_i(y_i,z_i)=|y_i-z_i|^{p}$, $\alpha=dC$, and $\varphi_i=F_{\mu_i}^{-1}$, $\psi_i=F_{\rho_i}^{-1}$ in the following Lemma.
\begin{alem}\label{lemprinc}
   Let $\alpha$ be a probability measure on $\R^{d_1}\times\R^{d_2}\times\hdots\times\R^{d_n}$ with respective marginals $\alpha_1,\hdots,\alpha_n$ on $\R^{d_1},\hdots,\R^{d_n}$ and $\varphi_i:\R^{d_i}\to\R^{d'_i}$, $\psi_i:\R^{d_i}\to\R^{d''_i}$ and $c_i:\R^{d'_i}\times\R^{d''_i}\to\R_+$ be measurable functions such that 
\begin{equation}\label{assumption_lemma}\forall i\in\{1,\hdots,n\},\;\inf_{\nu_i<^{\alpha_i\circ \varphi_i^{-1}}_{\alpha_i\circ \psi_i^{-1}}}\int_{\R^{d'_i}\times\R^{d''_i}}c_i(y_i,z_i)\nu_i(dy_i,dz_i)=\int_{\R^{d_i}}c_i(\varphi_i(x_i),\psi_i(x_i))\alpha_i(dx_i).
\end{equation}
Then setting $\varphi:x=(x_1,\hdots,x_n)\ni\R^{d_1}\times\hdots\times\R^{d_n}\mapsto (\varphi_1(x_1),\hdots,\varphi_n(x_n))\in\R^{d'_1+\hdots+d'_n}$ and $\psi:x\ni\R^{d_1}\times\hdots\times\R^{d_n}\mapsto (\psi_1(x_1),\hdots,\psi_n(x_n))\in\R^{d''_1+\hdots+d''_n}$, one has
\begin{align*}
   \inf&_{\nu<^{\alpha\circ \varphi^{-1}}_{\alpha\circ \psi^{-1}}}\int_{\R^{d'_1+\hdots+d'_n}\times\R^{d''_1+\hdots+d''_n}}\sum_{i=1}^n c_i(y_i,z_i)\nu(dy,dz)\\&=\int_{\R^{d_1+\hdots+d_n}}\sum_{i=1}^n c_i(\varphi_i(x_i),\psi_i(x_i))\alpha(dx)=\sum_{i=1}^n\int_{\R^{d_i}}c_i(\varphi_i(x_i),\psi_i(x_i))\alpha_i(dx_i).
\end{align*}
\end{alem}
\begin{adem}[of Lemma \ref{lemprinc}]
We give two alternative proofs of the Lemma. The first one is based on basic arguments. 
\begin{align*}
\int_{\R^{d_1+\hdots+d_n}}\sum_{i=1}^n &c_i(\varphi_i(x_i),\psi_i(x_i))\alpha(dx)\\
&\geq\inf_{\nu<^{\alpha\circ \varphi^{-1}}_{\alpha\circ \psi^{-1}}}\int_{\R^{d'_1+\hdots+d'_n}\times\R^{d''_1+\hdots+d''_n}}\sum_{i=1}^n c_i(y_i,z_i)\nu(dy,dz)\\
&\geq  \sum_{i=1}^n \inf_{\nu<^{\alpha\circ \varphi^{-1}}_{\alpha\circ \psi^{-1}}}\int_{\R^{d'_1+\hdots+d'_n}\times\R^{d''_1+\hdots+d''_n}}c_i(y_i,z_i)\nu(dy,dz)\\
&\geq \sum_{i=1}^n\inf_{\nu_i<^{\alpha_i\circ \varphi_i^{-1}}_{\alpha_i\circ \psi_i^{-1}}}\int_{\R^{d'_i}\times\R^{d''_i}}c_i(y_i,z_i)\nu_i(dy_i,dz_i)\\
&=\sum_{i=1}^n\int_{\R^{d_i}}c_i(\varphi_i(x_i),\psi_i(x_i))\alpha_i(dx_i)\\
&=\int_{\R^{d_1+\hdots+d_n}}\sum_{i=1}^n c_i(\varphi_i(x_i),\psi_i(x_i))\alpha(dx),
\end{align*}
where we used that 
\begin{itemize}
   \item the probability measure $\alpha\circ(\varphi^{-1},\psi^{-1})$ on $\R^{d'_1+\hdots+d'_n}\times\R^{d''_1+\hdots+d''_n}$ has respective marginals $\alpha\circ \varphi^{-1}$ and $\alpha\circ \psi^{-1}$ on $\R^{d'_1+\hdots+d'_n}$ and $\R^{d''_1+\hdots+d''_n}$, for the first inequality,
\item the infimum of a sum is greater than the sum of infima, for the second inequality,
\item the respective marginals of $\alpha\circ \varphi^{-1}$ and $\alpha\circ \psi^{-1}$ on $\R^{d_i'}$ and $\R^{d''_i}$ are $\alpha_i\circ \varphi_i^{-1}$ and $\alpha_i\circ \psi_i^{-1}$, for the third one,
\item and the hypotheses for the first equality.
\end{itemize}

The second proof is given to illustrate the theory of optimal transport. It requires to make the following additional assumption on the cost function $c_i$, for all $1\le i \le n$:
\begin{align}&\text{there exist measurable functions }g_i:\R^{d'_i}\rightarrow \R_+ \text{ and }h_i:\R^{d''_i}\rightarrow \R_+ \text{ such that }\nonumber \\
& c_i(y_i,z_i)\le g_i(y_i)+h_i(z_i), \ \int_{\R^{d_i}} g_i(\varphi_i(x_i)) \alpha_i(dx_i)< \infty ,   \ \int_{\R^{d_i}} h_i(\psi_i(x_i))\alpha_i(dx_i) <\infty.
\end{align}
Basically, this assumption ensures that $\int_{\R^{d'_i}\times\R^{d''_i}}c_i(y_i,z_i)\nu_i(dy_i,dz_i)< \infty$ for any coupling $\nu_i\hspace{-1mm}<^{\alpha_i\circ \varphi_i^{-1}}_{\alpha_i\circ \psi_i^{-1}}$.

We now introduce some definitions that are needed, and refer to the Section~3.3 of Rachev and  R\"uschendorf~\cite{raru} for a full introduction. Let $\bar{c}:\R^{d'}\times \R^{d''} \rightarrow \R $. A function $f:\R^{d'}\rightarrow \bar{\R}$ is {\it $\bar{c}$-convex } if there is a function $a:\R^{d''}\rightarrow  \bar{\R}$ such that 
$f(y)=\sup_{z\in \R^{d''}} \bar{c}(y,z)-a(z)$. For $y\in \R^{d'}$, we define the $\bar{c}$-subgradient:
$$\partial_{\bar{c}}f(y)=\{z \in \R^{d''} \ s.t. \  \forall \tilde{y} \in \text{dom} f,\ f(\tilde{y})-f(y) \ge \bar{c}(\tilde{y},z)-\bar{c}(y,z) \},$$
where $\text{dom} f=\{ y \in \R^{d'}, f(y)<\infty \}$. 

We are now in position to prove the result again. Let $X=(X_1,\dots,X_d)$ be a random variable with probability measure~$\alpha$. We define $Y_i=\varphi_i(X_i)$, $Z_i=\psi_i(X_i)$ and $\bar{c}_i=-c_i$. From~\eqref{assumption_lemma}, we know that $(Y_i,Z_i)$  is an optimal coupling that maximizes $\E[\bar{c}_i(Y_i,Z_i)]$. {From} Theorem 3.3.11 of Rachev and Rüschendorf~\cite{raru}, this implies the existence of a $\bar{c}_i$-convex function~$f_i:\R^{d'_i}\rightarrow \bar{\R}$ such that $Z_i\in \partial_{\bar{c}_i}f_i(Y_i)$. By definition of the  $\bar{c}_i$-convexity, there exists a function $a_i:\R^{d''_i}\rightarrow \bar{\R}$ such that $f_i(y_i)=\sup_{z_i \in \R^{d''_i}} \bar{c}(y_i,z_i)-a(z_i)$.
 We define for $y=(y_1,\dots,y_n)\in\R^{d_1'}\times\hdots\times\R^{d_n'}$ and $z=(z_1,\dots,z_n)\in\R^{d_1''}\times\hdots\times\R^{d_n''}$, 
$$f(y)=\sum_{i=1}^n f_i(y_i), \ a(z)= \sum_{i=1}^n a_i(z_i)\text{ and } \bar{c}(y,z)=\sum_{i=1}^n \bar{c}_i(y_i,z_i).$$
The function $f$ is $\bar{c}$-convex since we clearly have $f(y)=\sup_{z\in \R^{d''_1+\dots+d''_n}}\bar{c}(y,z)-a(z)$. Then, we have the straightforward inclusion:
\begin{align*}
\partial_{\bar{c}} f(y)&=\{(z_1,\dots,z_n)\in\R^{d_1''}\times\hdots\times\R^{d_n''},  \\
& \ \ \ \ \ s.t. \  \forall \tilde{y} \in \text{dom} f,\ \sum_{i=1}^n f_i(\tilde{y}_i)- f_i(y_i)\ge \sum_{i=1}^n \bar{c}_i(\tilde{y}_i,z_i)-\sum_{i=1}^n \bar{c}_i(y_i,z_i)\}\\
&\supset \partial_{\bar{c}_1}f_1(y_1) \times \dots \times \partial_{\bar{c}_n}f_n(y_n).
\end{align*}
This gives immediately  $Y\in \partial_{\bar{c}}f(X)$. Using again Theorem 3.3.11 of~\cite{raru}, we get that the coupling $(Y,Z)$ with law $\mu \diamond \rho$ is optimal in the sense that it maximizes $\E[\bar{c}(Y,Z)]$. \end{adem}

We now prove that the coupling $\mu \diamond \rho$ is in general no longer optimal when $q\not = p$. We first deal with the dimension $n=2$. Given two copulas $C_2$ and $C'_2$ on $[0,1]^2$, let $(U_1,U'_2)$ be distributed according $dC_2'$. Given $(U_1,U'_2)$, let $U_2$ be distributed according to the conditional distribution of the second coordinate given that the first one is equal to $U_1$ under $dC_2$ and $U'_1$ be distributed according to the conditional distribution of the first coordinate given that the second one in equal to $U'_2$ still under $dC_2$. This way the random variables $U_1$, $U_2$, $U_1'$ and $U_2'$ are uniformly distributed on $[0,1]$ and both the vectors $(U_1,U_2)$ and $(U'_1,U'_2)$ are distributed according to $dC_2$.  

For $\varepsilon\in [0,1]$, we consider
$$Y_\varepsilon=(U_1,\varepsilon U_2), \  Z_\varepsilon= (\varepsilon U_1,U_2) \text{ and } Z_\varepsilon'= (\varepsilon U'_1,U'_2).$$
We notice that the copula of $Y_\varepsilon$ and $Z_\varepsilon$ is~$C_2$ since the copula is preserved by coordinatewise increasing functions (see Nelsen~\cite{nelsen}, Theorem 2.4.3). Also, $Z_\varepsilon$ and $Z'_\varepsilon$ obviously have the same law. We will show that for a suitable choice of $C_2'$ and $\varepsilon>0$ small enough, we generally have $$\E[\|Y_\varepsilon-Z_\varepsilon\|_q^p]>\E[\|Y_\varepsilon-Z'_\varepsilon\|_q^p].$$
Denoting by $\mu^\varepsilon$ the law of $Y_\varepsilon$ and $\rho^\varepsilon$ the common law of $Z_\varepsilon$ and $Z'_\varepsilon$, this implies the desired conclusion since
\begin{align*}
  \int_{\R^2\times\R^2}\|x-y\|_q^p &\mu^\varepsilon \diamond \rho^\varepsilon(dx,dy)=\int_{[0,1]^2}\|(F_{\mu^\varepsilon_1}^{-1}(u_1),F_{\mu^\varepsilon_2}^{-1}(u_2))-(F_{\rho^\varepsilon_1}^{-1}(u_1),F_{\rho^\varepsilon_2}^{-1}(u_2))\|_q^pdC_2(u_1,u_2)\\
&=\E[\|Y_\varepsilon-Z_\varepsilon\|_q^p]>\E[\|Y_\varepsilon-Z'_\varepsilon\|_q^p]\geq {\cal W}_{p,q}^p(\mu,\rho).
 \end{align*} 
The coupling between $Y_\varepsilon$ and $Z_\varepsilon$ gives the score
$$\E[((1-\varepsilon)^qU_1^q+ (1-\varepsilon)^qU_2^q)^{p/q}] \underset{\varepsilon \rightarrow 0}{\rightarrow } \E[(U_1^q+ U_2^q)^{p/q}],$$ 
while the one between $Y_\varepsilon$ and $Z'_\varepsilon$ gives:
$$\E[(|U_1-\varepsilon U'_1|^q+ |U'_2-\varepsilon U_2|^q)^{p/q}] \underset{\varepsilon \rightarrow 0}{\rightarrow } \E[(U_1^q+ (U_2')^q)^{p/q}].$$
We now focus on the cost function $\bar{c}(u_1,u_2)=-(u_1^q+ u_2^q)^{p/q}$ for $u_1,u_2 \in (0,1)$. We have 
$$\partial_{u_1}\partial_{u_2}\bar{c}(u_1,u_2)=p(q-p)u_1^{q-1} u_2^{q-1}(u_1^q+ u_2^q)^{\frac{p}{q}-2}.$$
When $q<p$ (resp.~$q>p$) this is negative (resp.~positive) for any $u_1,u_2\in (0,1)$, i.e. $\bar{c}$ (resp. $-\bar{c}$) satisfies the so-called Monge condition. By Theorem~3.1.2 of Rachev and  R\"uschendorf~\cite{raru}, we get that $\E[\bar{c}(U_1,U'_2)]$ is maximal for $U'_2=1-U_1$ (resp. $U'_2=U_1$), i.e when $C_2'(u,v)=C^-_2(u,v)$ (resp. $C_2'(u,v)=C^+_2(u,v)$). Besides, since $\partial_{u_1}\partial_{u_2}\bar{c}(u_1,u_2)$ does not vanish, we have 
$$\E[\bar{c}(U_1,U_2)]< \E[\bar{c}(U_1,1-U_1)] \ (\text{resp. }\E[\bar{c}(U_1,U_2)]< \E[\bar{c}(U_1,U_1)]) $$
when $C_2 \not = C^-_2$ (resp.~$C_2 \not = C^+_2$). Taking $U'_2=1-U_1$ (resp. $U'_2=U_1$), we have in this case
$$\E[((1-\varepsilon)^qU_1^q+ (1-\varepsilon)^qU_2^q)^{p/q}]> \E[(|U_1-\varepsilon U'_1|^q+ |U'_2-\varepsilon U_2|^q)^{p/q}] $$
for $\varepsilon>0$ small enough. Notice that since both $Y_0$ and $Z_0$ have one constant coordinate, the range of the cumulative distribution function of this coordinate is $\{0,1\}$ and the probability measures $\mu_0$ and $\rho_0$ share any two dimensional copula and in particular $C^-_2$ (resp. $C^+_2$). That is why one has to choose $\varepsilon>0$.

This two dimensional example  can be easily extended to dimension~$n\ge 3$. Let $C_n$ now denote a $n$-dimensional copula, $C_2(u_1,u_2)=C_n(u_1,u_2,1,\dots,1)$ and $C'_2$ be another two-dimensional copula. We define $(U_1,U_2,U'_1,U'_2)$ as above. Then, we choose
$(U_3,\dots,U_n)$ (resp.~$(U'_3,\dots,U'_n)$)  distributed according to the conditional law of the $n-2$ last coordinates given that the two first are equal to $(U_1,U_2)$ (resp.~$(U'_1,U'_2)$) under $dC_n$. Last, we define $$Y_\varepsilon=(U_1,\varepsilon U_2,\varepsilon U_3,\dots,\varepsilon U_n), \ Z_\varepsilon=(\varepsilon U_1,U_2,\varepsilon U_3,\dots,\varepsilon U_n) \text{ and } Z'_\varepsilon=(\varepsilon U'_1,U'_2,\varepsilon U'_3,\dots,\varepsilon U'_n).$$
We still have on the one hand that $Y_\varepsilon$ and $Z_\varepsilon$ share the same copula and on the other hand that $Z_\varepsilon$ and $Z'_\varepsilon$ have the same law. Moreover, 
$$\E[\|Y_\varepsilon-Z_\varepsilon\|_q^p]\underset{\varepsilon \rightarrow 0}{\rightarrow } \E[(U_1^q+ U_2^q)^{p/q}] \text{ and } \E[\|Y_\varepsilon-Z'_\varepsilon\|_q^p]\underset{\varepsilon \rightarrow 0}{\rightarrow } \E[(U_1^q+ (U_2')^q)^{p/q}]. $$
Taking $\varepsilon>0$ small enough, we get that  $\E[\|Y_\varepsilon-Z_\varepsilon\|_q^p] > \E[\|Y_\varepsilon-Z'_\varepsilon\|_q^p]\geq {\cal W}_{p,q}^p(\mu,\rho)$ if, for some $u_1,u_2 \in [0,1]$, $C_n(u_1,u_2,1,\dots,1)>C_2^-(u_1,u_2)$ when $q<p$ (resp. $C_n(u_1,u_2,1,\dots,1)<C_2^+(u_1,u_2)$ when $q>p$).  
If there exist $i<j$ such that
$$C_n(1,\dots,1,u_i,1,\dots,1,u_j,1,\dots,1)>C_2^-(u_i,u_j)  \ (\text{resp. } < C_2^+(u_i,u_j)\ ),  $$
then one may repeat the above reasoning with the $i$-th and $j$-th coordinates replacing the first and second ones. Then the coupling $\nu=\mu \diamond \rho$ is not optimal for $\varepsilon>0$ small enough. 
 For $n\ge 3$, there is no copula such that $C_n(1,\dots,1,u_i,1,\dots,1,u_j,1,\dots,1)  = C_2^-(u_i,u_j)$  for any $ i<j$ , $u_i,u_j\in[0,1]$, since this would imply that $U_2=1-U_1=U_3$ and $U_3=1-U_2$. Also, the only one copula satisfying $C_n(1,\dots,1,u_i,1,\dots,1,u_j,1,\dots,1)  =  C_2^+(u_i,u_j)$ for any $ i<j$, $u_i,u_j\in[0,1]$,  is $C_n^+$ since the former condition implies $U_i=U_j$ for any $i<j$.

\end{document}